\documentclass{article}
\usepackage{amsfonts}
\usepackage{amssymb, amsmath, mathrsfs, amsfonts, cite}
\begin{document}

 \title{\bf Twisted Angles for Central Configurations Formed By Two Twisted Regular Polygons\thanks{Supported partially by NSF of China}}
  \author{Yu Xiang and Zhang Shiqing\\
  {\small Yantze Center of Mathematics, Sichuan University, Chengdu 610064,P.R.China}\\
  }

 \date{}
 \maketitle
(Dedicated to the Memory of Professor Chern S.S. on the Occasion of His 100th Birthday)\\\\
\noindent{\small {\bf Abstract:} In this paper, we study the necessary conditions and sufficient conditions for the twisted angles of the central configurations formed by two twisted regular polygons, specially, we prove that for the 2N-body problems, the twisted angles must be
$\theta=0~\mbox{or}~\theta=\pi/N$.

 \vskip 2mm


 \vskip 2mm
 \noindent {\bf 2000AMS Mathematical Subject Classification:} 70F10,70F15.\\

\section{Introduction and Main Results}

The problem on the numbers for central configurations is so important that S.Smale ([7])took it as one of the most important 18 mathematical problems for the 21'st century. Central configurations are important in the Newtonian N-body problems, for example, it's well known that finding the relative equilibrium solutions of the classical N-body problem and planar central configurations is equivalent.Central configurations also play other important roles, there were a lot of works on the existence and multiplicity and the shapes of central configurations ([6]). But finding concrete central configurations is very difficult, here we consider some particular situations: the central configurations formed by two twisted regular polygons, we will explain it more clearly in the following.

{\bf Definition ([6,8]):} A configuration $q=(q_1,\ldots,q_n)\in \textit{X}\setminus\Delta$ is called a central configuration if there exists a constant $\lambda\in \textbf{R}$ such that
\begin{equation}
\sum_{j=1,j\neq k}^n \frac{m_jm_k}{|q_j-q_k|^3}(q_j-q_k)=-\lambda m_kq_k,1\leq k\leq n
\end{equation}

The value of $\lambda$ in (1.1) is uniquely determined by
\begin{equation}
\lambda=\frac{U(q)}{I(q)}
\end{equation}

Where
\begin{eqnarray}
\textit{X}=\left\{q=(q_1,q_2,\ldots,q_n)\in \textbf{R}^{3n}:\sum_{i=1}^n m_iq_i=0\right\}\\
\Delta=\left\{q:q_j=q_k~\mbox{for~some}~j\neq k\right\}~~~~~~~~~~~\\
U(q)=\sum_{1\leq j<k\leq n} \frac{m_jm_k}{|q_j-q_k|}~~~~~~~~~~~~~~~~~\\
I(q)=\sum_{1\leq j\leq n} m_j|q_j|^2~~~~~~~~~~~~~~~~~~~
\end{eqnarray}

Consider the central configurations in $\textbf{R}^3$ formed by two twisted regular N-gons $(N\geq 2)$ with distance $h\geq0$.It is assumed that the lower layer regular N-gons lies in horizontal plane, and the upper regular N-gons parallels to the lower one and z-axis passes through both centers of two regular N-gons. Suppose that the lower layer particles have masses $m_1,m_2,\ldots,m_N$ and the upper layer particles have masses $\tilde{m}_1,\tilde{m}_2,\ldots,\tilde{m}_N$ respectively. More precisely, let $\rho_k$ be the k-th complex root of unity, i.e.,
\begin{equation}
\rho_k=e^{i\theta_k}
\end{equation}
And we let
\begin{equation}
\tilde{\rho}_k=a\rho_k\cdot e^{i\theta}
\end{equation}
Where $a>0,i=\sqrt{-1},\theta_k=\frac{2k\pi}{N}(k=1,\ldots,N),0\leq \theta\leq 2\pi,\theta$ is called twisted angle.

It is assumed that $m_k(k=1,\ldots,N)$ locates at the vertex $q_k$ of the lower layer regular N-gons;$\tilde{m}_k(k=1,\ldots,N)$ 	locate at the vertex of the upper layer regular N-gons:
\begin{eqnarray}
q_k=(\rho_k,0)\\
\tilde{q}_k=(\tilde{\rho}_k,0)
\end{eqnarray}
Where $h\geq 0$ is the distance between the two layers.
Then the center of masses is
\begin{equation}
z_0=\frac{\sum_j (m_jq_j+\tilde{m}_j\tilde{q}_j)}{M}
\end{equation}
Where
\begin{equation}
M=\sum_j (m_j+\tilde{m}_j)
\end{equation}
Let
\begin{eqnarray}
P_k=q_k-z_0,P=(P_1,\ldots,P_N)\\
\tilde{P}_k=\tilde{q}_k-z_0,\tilde{P}=(\tilde{P}_1,\ldots,\tilde{P}_N)
\end{eqnarray}
If $P_1,\ldots,P_N;\tilde{P}_1,\ldots,\tilde{P}_N$ form a central configuration,then $\exists~\lambda \in \textbf{R}^+$,such that\\
\begin{eqnarray}
\sum_{j=1,j\neq k}^N \frac{m_j}{|P_k-P_j|^3}(P_k-P_j)+\sum_{j=1}^N \frac{\tilde{m}_j}{|P_k-\tilde{P}_j|^3}(P_k-\tilde{P}_j)=\lambda P_k,1\leq k\leq N\\
\sum_{j=1,j\neq k}^N \frac{\tilde{m}_j}{|\tilde{P}_k-\tilde{P}_j|^3}(\tilde{P}_k-\tilde{P}_j)+\sum_{j=1}^N \frac{m_j}{|\tilde{P}_k-P_j|^3}(\tilde{P}_k-P_j)=\lambda \tilde{P}_k,1\leq k\leq N\\
\lambda=\frac{U(P,\tilde{P})}{I(P,\tilde{P})}~~~~~~~~~~~~~~~~~~~~~~~~~~~~~~~~~~~~~~~~~~~~
\end{eqnarray}\\
In the following, we only consider the case of $m_1=\cdots=m_N=m$ and $\tilde{m}_1=\cdots=\tilde{m}_N=bm$.
Then\\
\begin{equation}
z_0=\sum_j (m_jq_j+\tilde{m}_j\tilde{q}_j)/M=\left(0,0,\frac{bh}{1+b}\right)
\end{equation}\\
And (1.15),(1.16)and(1.17)are equivalent to\\
\begin{equation}
\sum_{1\leq j\leq N-1} \frac{(1-\rho_j,0)}{|1-\rho_j|^3}+b\sum_{1\leq j\leq N} \frac{(1-\tilde{\rho}_j,-h)}{(|1-\tilde{\rho}_j|^2+h^2)^\frac{3}{2}}=\frac{\lambda}{m}\left(1,0,-\frac{bh}{1+b}\right)
\end{equation}
\begin{eqnarray}
\frac{b}{a^3}\sum_{1\leq j\leq N-1} \frac{(ae^{i\theta}-\tilde{\rho}_j,0)}{|1-\rho_j|^3}+\sum_{1\leq j\leq N} \frac{(ae^{i\theta}-\rho_j,h)}{(|ae^{i\theta}-\rho_j|^2+h^2)^\frac{3}{2}}=\frac{\lambda}{m}\left(ae^{i\theta},\frac{h}{1+b}\right)\\\nonumber\\
\frac{\lambda}{m}=\frac{\left(1+\frac{b^2}{a}\right)\sum_{1\leq j<k\leq N} \frac{1}{|\rho_j-\rho_k|}+bN\sum_{1\leq j\leq N} \frac{1}{(|1-\tilde{\rho}_j|^2+h^2)^{1/2}}}{N\left(1+ba^2+\frac{bh^2}{1+b}\right)}~~~~~~~~~~~~
\end{eqnarray}\\
In the following,let $\mu=\frac{\lambda}{m}$,and we only consider the case of $0\leq \theta\leq \frac{2\pi}{N}$ or $-\frac{\pi}{N}\leq \theta\leq \frac{\pi}{N}$ because of the symmetry.\\
When $\theta=0$,R. Moeckel and C. Simo ([4]) proved the following results:\\

{\bf Theorem 1.1}(R. Moeckel and C.Simo). When $h=0,\theta=0$,for every mass ratio $b$,there are exactly two planar central configurations consisting of two nested regular N-gons. For one of these, the ratio $a$ of the sizes of the two polygons is less than 1, and for the other it is greater than 1. However, for $N\geq 473$ there is a constant $b_0(N)<1$ such that for $b<b_0$ and $b>\frac{1}{b_0}$,the central configuration with the smaller masses on the inner polygon is a repeller.\\

{\bf Theorem 1.2}(R. Moeckel and C.Simo). When $h>0,\theta=0$,if $N<473$,there is a unique pair of spatial central configurations of parallel regular N-gons. If $N\geq 473$,here are no such central configurations for $b<b_0(N)$.At $b=b_0$ a unique pair bifurcates from the planar central configuration with the smaller masses on the inner polygon. This remains the unique pair of spatial central configurations until $b=\frac{1}{b_0}$,where a similar bifurcation occurs in reverse, so that for $b>\frac{1}{b_0}$,only the planar central configurations remain.\\

	Zhang-Zhou ([10]) studied the necessary and sufficient conditions for the masses of the central configurations consisting of two planar twisted regular N-gons. They proved the following Theorem:\\

{\bf Theorem1.3} If the central configuration is formed by two twisted regular N-gons $(N\geq2)$ with distance $h=0$, then\\
\begin{eqnarray}
\lambda\frac{N}{M}=\frac{1}{1+b}\left(\sum_{1\leq j\leq N-1} \frac{1-\rho_j}{|1-\rho_j|^3}+\sum_{1\leq j\leq N} \frac{b\left(1-a\rho_je^{i\theta}\right)}{|1-a\rho_je^{i\theta}|^3}\right)~~\\\nonumber\\
\lambda\frac{N}{M}=\frac{e^{-i\theta}}{a(1+b)}\left(\sum_{1\leq j\leq N-1} \frac{b(1-\rho_j)e^{i\theta}}{a^2|1-\rho_j|^3}+\sum_{1\leq j\leq N} \frac{ae^{i\theta}-\rho_j}{|ae^{i\theta}-\rho_j|^3}\right)
\end{eqnarray}\\
From their Theorem, a series of conclusions are derived, especially they have the following Corollaries :\\

{\bf Corollary1.4}(MacMillan- Bartky [3])~~For $N=2$ and $\theta=\pi/2$,then $b=1$ if and only if $a=1$.\\

{\bf Corollary1.5}(Perko-Walter [5])~~For $N\geq2$,$a=1$ and $\theta=\pi/N$,if (1.15), (1.16) and (1.17) hold, then $b=1$.\\

When $h>0$ 	Xie-Zhang-Zhou([9]) proved:\\

{\bf Theorem1.6}(Xie-Zhang-Zhou)~~ The configuration formed by two twisted regular N-gons $(N\geq2)$ with distance $h>0$ is a central configuration if and only if the parameters $a,b,h,\theta$ satisfy the following relationships:\\
\begin{equation}
\lambda\frac{N}{M}=\frac{1}{1+b}\left(\sum_{1\leq j\leq N-1} \frac{1-\rho_j}{|1-\rho_j|^3}+\sum_{1\leq j\leq N} \frac{b(1-ae^{i\theta}\rho_j)}{\left(|1-ae^{i\theta}\rho_j|^2+h^2\right)^{3/2}}\right)
\end{equation}
\begin{eqnarray}
\lambda\frac{N}{M}=\frac{e^{-i\theta}}{a(1+b)}\left(\sum_{1\leq j\leq N-1} \frac{b(1-\rho_j)e^{i\theta}}{a^2|1-\rho_j|^3}+\sum_{1\leq j\leq N} \frac{ae^{i\theta}-\rho_j}{\left(|ae^{i\theta}-\rho_j|^2+h^2\right)^{3/2}}\right)\\\nonumber\\
\lambda\frac{N}{M}=\sum_{1\leq j\leq N} \frac{1}{\left(|1-ae^{i\theta}\rho_j|^2+h^2\right)^{3/2}}~~~~~~~~~~~~~~~~~~~~~~~~~~
\end{eqnarray}\\
Zhang-Zhu ([11]) proved the following result for a special case:\\

{\bf Theorem 1.7}(Zhang-Zhu) When $h>0$,if $a=1,b=1$,and $\theta=\pi/N$,then for every $N$,	there exists a unique central configuration. Particularly, they proved that
\begin{equation}
\sum_{1\leq j\leq N} \frac{\cos(\theta_j+\pi/N)+1}{\left(2-2\cos(\theta_j+\pi/N)\right)^{3/2}}-\sum_{1\leq j\leq N-1} \frac{1-\rho_j}{|1-\rho_j|^3}>0
\end{equation}\\

In the following, we let
\begin{equation}
A=\sum_{1\leq j\leq N-1} \frac{1-\rho_j}{|1-\rho_j|^3}
\end{equation}\\
In this paper we will prove the following main results:\\

{\bf Theorem 1.8}~~~If the central configuration is formed by two twisted regular N-gons $(N\geq2)$ with distance $h\geq0$, then only $\theta=0$ or $\theta=\pi/N$.Specially,if $a=1$ and $h=0$,i.e.,two nested regular N-gons are on the same unit circle, then only $\theta=\pi/N$.\\

{\bf Corollary 1.9}~~~For $N\geq2,h=0$,if $a=1$, then $b=1$ and $\theta=\pi/N$,i.e., there is exactly one central configuration formed by two nested regular N-gons on the same unit circle, which is the regular 2N-gons.\\

{\bf Corollary 1.10}~~~The configuration formed by two twisted regular N-gons $(N\geq2)$ with diatance $h\geq0$ is a central configuration if and only if the parameters $a,b,h$ satisfy the following relationships:\\\\
i.~~~When $h=0$ and $a\neq1$\\
\begin{equation}b\left[\sum_{1\leq j\leq N} \frac{1-a\cos(\theta_j)}{(1+a^2-2a\cos(\theta_j))^{3/2}}-\frac{A}{a^3}\right]=\sum_{1\leq j\leq N} \frac{1-a^{-1}\cos(\theta_j)}{(1+a^2-2a\cos(\theta_j))^{3/2}}-A \nonumber
\end{equation}\\
or\\
\begin{equation}b\left[\sum_{1\leq j\leq N} \frac{1-a\cos(\theta_j+\frac{\pi}{N})}{(1+a^2-2a\cos(\theta_j+\frac{\pi}{N}))^{3/2}}-\frac{A}{a^3}\right]=\sum_{1\leq j\leq N} \frac{1-a^{-1}\cos(\theta_j+\frac{\pi}{N})}{(1+a^2-2a\cos(\theta_j+\frac{\pi}{N}))^{3/2}}-A \nonumber
\end{equation}\\
ii.When $h>0$\\
\begin{eqnarray}
\left\{
\begin{array}{c}
ba\sum_{1\leq j\leq N} \frac{\cos(\theta_j)}{(1+a^2-2a\cos(\theta_j)+h^2)^{3/2}}=A-\sum_{1\leq j\leq N} \frac{1}{(1+a^2-2a\cos(\theta_j)+h^2)^{3/2}}\nonumber\\\nonumber\\
ba\left(\frac{A}{a^3}-\sum_{1\leq j\leq N} \frac{1}{(1+a^2-2a\cos(\theta_j)+h^2)^{3/2}}\right)=\sum_{1\leq j\leq N} \frac{\cos(\theta_j)}{(1+a^2-2a\cos(\theta_j)+h^2)^{3/2}}\nonumber
\end{array}
\right.
\end{eqnarray}\\
or\\
\begin{eqnarray}
\left\{
\begin{array}{c}
ba\sum_{1\leq j\leq N} \frac{\cos(\theta_j+\frac{\pi}{N})}{(1+a^2-2a\cos(\theta_j+\frac{\pi}{N})+h^2)^{3/2}}=A-\sum_{1\leq j\leq N} \frac{1}{(1+a^2-2a\cos(\theta_j+\frac{\pi}{N})+h^2)^{3/2}}\nonumber\\\nonumber\\
ba\left(\frac{A}{a^3}-\sum_{1\leq j\leq N} \frac{1}{(1+a^2-2a\cos(\theta_j+\frac{\pi}{N})+h^2)^{3/2}}\right)=\sum_{1\leq j\leq N} \frac{\cos(\theta_j)}{(1+a^2-2a\cos(\theta_j+\frac{\pi}{N})+h^2)^{3/2}}\nonumber
\end{array}
\right.
\end{eqnarray}\\

{\bf Corollary 1.11}~~~For $N\geq2,h>0,a=1$,if the configuration formed by two twisted regular N-gons $(N\geq2)$ with distance $h\geq0$ is a central configuration, then $b=1,\theta=0~\mbox{or}~\pi/N$,and there exists a unique $h$ for each $\theta$.	In other words, there are exactly two spatial central configurations formed by parallel regular N-gons which have the same sizes.
\setcounter{equation}{0}

\section{Some Lemmas}

First of all, let's establish some identical equations to simplify the problem.

\begin{flalign}
&~~~~~\mbox{{\bf Lemma 2.1}}~~\frac{1}{N}\sum_{1\leq j<k\leq N} \frac{1}{|\rho_j-\rho_k|}=\sum_{1\leq j\leq N-1} \frac{1-\rho_j}{|1-\rho_j|^3}=\frac{1}{4}\sum_{1\leq j\leq N-1} \csc\left(\frac{\pi j}{N}\right)&
\end{flalign}\\

{\bf Proof:} It's easy to know
\begin{equation}\sum_{1\leq j\leq N-1} \frac{1-\rho_j}{|1-\rho_j|^3}=\frac{1}{4}\sum_{1\leq j\leq N-1} \csc\left(\frac{\pi j}{N}\right)>0.\nonumber
\end{equation}
We only need to prove
\begin{equation}
\sum_{1\leq j<k\leq N} \frac{1}{|\rho_j-\rho_k|}=\frac{N}{4}\sum_{1\leq j\leq N-1} \csc\left(\frac{\pi j}{N}\right).\nonumber\end{equation}
In fact, we have
\begin{equation}2\sum_{1\leq j<k\leq N} \frac{1}{|\rho_j-\rho_k|}=\sum_{j\neq k} \frac{1}{|\rho_j-\rho_k|}=N\sum_{1\leq j\leq N-1} \frac{1}{\rho_j-1}=\frac{N}{2}\sum_{1\leq j\leq N-1} \csc\left(\frac{\pi j}{N}\right).\nonumber~~~~~~~~~~~~~~~~~~~~~~~~~~~~~~~~~~~~\Box
\end{equation}\\
Then we have
\begin{equation}
A=\frac{1}{N}\sum_{1\leq j<k\leq N} \frac{1}{|\rho_j-\rho_k|}=\sum_{1\leq j\leq N-1} \frac{1-\rho_j}{|1-\rho_j|^3}=\frac{1}{4}\sum_{1\leq j\leq N-1} \csc\left(\frac{\pi j}{N}\right).
\end{equation}

{\bf Lemma 2.2}~~~ We have\\
\begin{eqnarray}
\sum_{1\leq j\leq N} \frac{cos(\theta_j+\theta)}{(1+a^2-2acos(\theta_j+\theta)+h^2)^{3/2}}=\sum_{1\leq j\leq N} \frac{cos(\theta_j-\theta)}{(1+a^2-2acos(\theta_j-\theta)+h^2)^{3/2}}\\\nonumber\\
\sum_{1\leq j\leq N} \frac{1}{(1+a^2-2acos(\theta_j+\theta)+h^2)^{3/2}}=\sum_{1\leq j\leq N} \frac{1}{(1+a^2-2acos(\theta_j-\theta)+h^2)^{3/2}}\\\nonumber\\
\sum_{1\leq j\leq N} \frac{cos(\theta_j+\theta)}{(1+a^2-2acos(\theta_j+\theta)+h^2)^{1/2}}=\sum_{1\leq j\leq N} \frac{1}{(1+a^2-2acos(\theta_j-\theta)+h^2)^{1/2}}\\\nonumber\\
\sum_{1\leq j\leq N} \frac{sin(\theta_j+\theta)}{(1+a^2-2acos(\theta_j+\theta)+h^2)^{3/2}}=-\sum_{1\leq j\leq N} \frac{sin(\theta_j-\theta)}{(1+a^2-2acos(\theta_j-\theta)+h^2)^{3/2}}
\end{eqnarray}

{\bf Proof:} For (2.3), we notice that
\begin{eqnarray}
\sum_{1\leq j\leq N} \frac{cos(\theta_j+\theta)}{(1+a^2-2acos(\theta_j+\theta)+h^2)^{3/2}}~~~~~~~~~~~~~~~~\nonumber\\
=\sum_{1\leq j\leq N} \frac{cos(2\pi-\theta_j-\theta)}{(1+a^2-2acos(2\pi-\theta_j-\theta)+h^2)^{3/2}}~~~~~~\nonumber\\
=\sum_{1\leq j\leq N} \frac{cos(\theta_{N-j}-\theta)}{(1+a^2-2acos(\theta_{N-j}-\theta)+h^2)^{3/2}}~~~~~~~~~\nonumber\\
=\sum_{1\leq k\leq N} \frac{cos(\theta_k-\theta)}{(1+a^2-2acos(\theta_k-\theta)+h^2)^{3/2}}~~~~~~~~~~~~\nonumber\\
=\sum_{1\leq k\leq N}
\frac{cos(\theta_k-\theta)}{(1+a^2-2acos(\theta_k-\theta)+h^2)^{3/2}}~~~~~~~~~~~~\nonumber
\end{eqnarray}\\
Similarly, we can get (2.4), (2.5) and (2.6).~~~~~~~~~~~~~~~~~~~~~~~~~~~~~~~~~~~~~~~~~~~~~~~~~~~~~~~~~~~~~~~~~~~~~~~~~~~~~~~~~~~~~$\Box$\\\\
From (1.19), (1.20) and (1.21), and using (2.1)-(2.6) we can get five equivalent equations .\\

{\bf Lemma 2.3} (1.19), (1.20) and (1.21) can be simplified into the following equations:\\
\begin{eqnarray}
A+b\sum_{1\leq j\leq N} \frac{1-a cos(\theta_j+\theta)}{(1+a^2-2acos(\theta_j+\theta)+h^2)^{3/2}}=\mu\\\nonumber\\
\frac{b}{a^3}A+\sum_{1\leq j\leq N} \frac{1-a^{-1}cos(\theta_j+\theta)}{(1+a^2-2acos(\theta_j+\theta)+h^2)^{3/2}}=\mu\\\nonumber\\
\sum_{1\leq j\leq N} \frac{sin(\theta_j+\theta)}{(1+a^2-2acos(\theta_j+\theta)+h^2)^{3/2}}=0~~~~~\\\nonumber\\
h\sum_{1\leq j\leq N} \frac{1}{(1+a^2-2acos(\theta_j+\theta)+h^2)^{3/2}}=\frac{\mu h}{1+b}\\\nonumber\\
\mu=\frac{\left(1+\frac{b^2}{a}\right)A+b\sum_{1\leq j\leq N} \frac{1+a^2+h^2-2acos(\theta+\theta)}{(1+a^2-2acos(\theta_j+\theta)+h^2)^{3/2}}}{\left(1+ba^2+\frac{bh^2}{1+b}\right)}
\end{eqnarray}\\

{\bf Proof:}We write every vector in the equations of (1.19) and (1.20) into the forms of components, and apply (2.1)-(2.6), then (2.7)-(2.10) is obvious. (2.11) is also clear from (2.1)-(2.6) and (1.21).~~~~~~~~~~~~~~~~~~~~~~~~~~~~~~~~~~~~~~~~~~~~~~~$\Box$\\\\
When $h=0$,we have the next Lemma:\\

{\bf Lemma 2.4} (1.19), (1.20) and (1.21) can be simplified into the following equations:\\
\begin{eqnarray}
A+b\sum_{1\leq j\leq N} \frac{1-acos(\theta_j+\theta)}{(1+a^2-2acos(\theta_j+\theta))^{3/2}}=\mu\\\nonumber\\
\frac{b}{a^3}A+\sum_{1\leq j\leq N} \frac{1-a^{-1}cos(\theta_j+\theta)}{(1+a^2-2acos(\theta_j+\theta))^{3/2}}=\mu
\end{eqnarray}
\begin{equation}
\sum_{1\leq j\leq N} \frac{sin(\theta_j+\theta)}{(1+a^2-2acos(\theta_j+\theta))^{3/2}}=0
\end{equation}\\

{\bf Proof:} We only indicate that (2.11) can be gotten from (2.7) and (2.8) when $h=0$.$~~~~~~~~~~~~~~~~~~~~~~~~~~~~~~~\Box$\\\\
When $h>0$,we have\\

{\bf Lemma 2.5} For $h>0$,	(1.19), (1.20) and (1.21) can be simplified into the following equations:\\
\begin{eqnarray}
A-ab\sum_{1\leq j\leq N} \frac{cos(\theta_j+\theta)}{(1+a^2-2acos(\theta_j+\theta)+h^2)^{3/2}}=\frac{\mu}{1+b}~~\\\nonumber\\
\frac{b}{a^3}A-a^{-1}\sum_{1\leq j\leq N} \frac{cos(\theta_j+\theta)}{(1+a^2-2acos(\theta_j+\theta)+h^2)^{3/2}}=\frac{b\mu}{1+b}\\\nonumber\\
\sum_{1\leq j\leq N} \frac{sin(\theta_j+\theta)}{(1+a^2-2acos(\theta_j+\theta)+h^2)^{3/2}}=0~~~~~~~~~~~\\\nonumber\\
\sum_{1\leq j\leq N} \frac{1}{(1+a^2-2acos(\theta_j+\theta)+h^2)^{3/2}}=\frac{\mu}{1+b}~~~~~~~~~~
\end{eqnarray}\\

{\bf Proof:} It is obvious by (2.7)-(2.10). We only indicate that (2.11) can be gotten from (2.7) and (2.8) and (2.18). $~~~~~~~~~~~~~~~~~~~~~~~~~~~~~~~~~~~~~~~~~~~~~~~~~~~~~~~~~~~~~~~~~~~~~~~~~~~~~~~~~~~~~~~~~~~~~~~~~~~~~~~~~~~~~~~~~~~~~~~~~~~~~~~~~~~~~~~~~~~~~~~~~~~~\Box$\\

{\bf Lemma 2.6([8])} For $n\geq 3,m_1=\cdots=m_n$,if $m_1,\cdots,m_n$ locate at vertices of a regular polygon ,then they form a central configuration.\\

{\bf Lemma 2.7} Let \begin{equation}g_n(x)=\sum_{1\leq j\leq N} \frac{sin(\theta_j+\theta)}{(1+a^2-2acos(\theta_j+\theta)+x)^{(2n+3)/2}},\nonumber\end{equation}
where $\theta\in\left(0,\frac{\pi}{N}\right)$ and $a>0,n \in \mathbb{N},x\geq0$.If $g_n(x)>0$ in $\{x:x\geq0\}$ for some $n\geq1$,then $g_j(x)>0$ in $\{x:x\geq0\}$ for $0\leq j\leq n-1$.\\

{\bf Proof :} First of all, it is easy to know that $g_m(x)\rightarrow0$ when $x\rightarrow\infty$ for any $m \in \mathbb{N}$.Since \begin{equation}
g'_{n-1}(x)=\left(-\frac{2n+1}{2}\right)\sum_{1\leq j\leq N}  \frac{sin(\theta_j+\theta)}{(1+a^2-2acos(\theta_j+\theta)+x)^{(2n+3)/2}}=\left(-\frac{2n+1}{2}\right)g_n(x)<0.
\nonumber\end{equation}
Thus $g_{n-1}(x)>0$ in $\{x:x\geq0\}$.
Similarly, we can get $g_{n-2}(x)>0$ in $\{x:x\geq0\}$,$\cdots$,$g_0(x)>0$ in \\ $\{x:x\geq0\}$.$~~~~~~~~~~~~~~~~~~~~~~~~~~~~~~~~~~~~~~~~~~~~~~~~~~~~~~~~~~~~~~~~~~~~~~~~~~~~~~~~~~~~~~~~~~~~~~~~~~~~~~~~~~~~~~~~~~~~~~~~~~~~~~~~~~~~~\Box$\\

{\bf Lemma 2.8} Let $a_j>0,1\leq j\leq k,A_1\geq\cdots\geq A_k\geq0$.Then $\lim_{n \rightarrow \infty} \left(\sum_{1\leq j\leq k} a_jA_j^n\right)^{\frac{1}{n}}=A_1$\\

{\bf Proof :} Let
$$B\geq a_j\geq b>0~\mbox{for}~1\leq j\leq k.$$
Then
$$b^{\frac{1}{n}}A_1\leq \left(\sum_{1\leq j\leq k} a_jA_j^n\right)^{\frac{1}{n}}\leq (kB)^{\frac{1}{n}}A_1.$$
So $\lim_{n\rightarrow\infty} \left(\sum_{1\leq j\leq k} a_jA_j^n\right)^{\frac{1}{n}}=A_1$ $~~~~~~~~~~~~~~~~~~~~~~~~~~~~~~~~~~~~~~~~~~~~~~~~~~~~~~~~~~~~~~~~~~~~~~~~~~~~~~~~~~~~~~~~~~~~~~~~~~\Box$\\

{\bf Lemma 2.9} Let $a,b>0$.Then $a-b$ and $a^{\frac{1}{n}}-b^{\frac{1}{n}}$ are positive or negative at the same time.\\\\
The proof of Lemma 2.9 is obviously.\\\\
{\bf Lemma 2.10} Let $f(\theta)=\sum_{1\leq j\leq N} \frac{\sin(\theta_j+\theta)}{|1+a^2-2a\cos(\theta_j+\theta)+h^2|^{3/2}}$,then $f(\frac{\pi}{N})=0$,$f(-\theta)=-f(\theta)$ and $f(\theta+\frac{2\pi}{N})=f(\theta)$.\\

{\bf Proof:}
\begin{eqnarray}
f(\frac{\pi}{N}) & = & \sum_{1\leq j\leq N} \frac{\sin((2j+1)\pi/N)}{|1+a^2-2a\cos((2j+1)\pi/N)+h^2|^{3/2}}\nonumber\\
& = & -\sum_{1\leq j\leq N} \frac{\sin((2N-2j-1)\pi/N)}{|1+a^2-2a\cos((2n-2j-1)\pi/N)+h^2|^{3/2}}\nonumber\\
& = & -\sum_{-1\leq k\leq N-2} \frac{\sin((2k+1)\pi/N)}{|1+a^2-2a\cos((2k+1)\pi/N)+h^2|^{3/2}}\nonumber\\
& = & -\sum_{1\leq k\leq N}
\frac{\sin((2k+1)\pi/N)}{|1+a^2-2a\cos((2k+1)\pi/N)+h^2|^{3/2}}\nonumber\\
& = & -f(\frac{\pi}{N})\nonumber
\end{eqnarray}
\begin{flalign}
&\mbox{Thus},f( \frac{\pi}{N} )=0.&
\end{flalign}
From (2.6) we have $f(-\theta)=-f(\theta)$,and $f(\theta+\frac{2\pi}{N})=f(\theta)$ is obvious by the definition of $f(\theta)$.$~~~~~~~~~~~~~~~~~~~~~~~~~~~~~~~~~~~~~~~~~~~\Box$\\

{\bf Lemma 2.11} We have $f(\theta)=\sum_{1\leq j\leq N} \frac{sin(\theta_j+\theta)}{|1+a^2-2a\cos(\theta_j+\theta)+h^2|^{3/2}}>0$ for any $a>0$ and $h\geq0$ when $\theta \in (0,\frac{\pi}{N})$.\\

{\bf Proof :}It is easy to know that we only need to prove
$$g_0(x)=\sum_{1\leq j\leq N} \frac{sin(\theta_j+\theta)}{(1+a^2-2a\cos(\theta_j+\theta)+x)^{3/2}}>0~\mbox{in}~ \{x:x\geq0\}$$ for any $\theta \in (0,\frac{\pi}{N})$ and $a>0$.But by Lemma 2.7 we know that we only need to prove that
$$g_n(x)=\sum_{1\leq j\leq N} \frac{sin(\theta_j+\theta)}{(1+a^2-2a\cos(\theta_j+\theta)+x)^{(2n+3)/2}}>0~\mbox{in}~\{x:x\geq0\}$$
for some sufficiently large $n \in \mathbb{N}$.Since
\begin{eqnarray}
g_n(x) & = & \sum_{0\leq j\leq N-1} \frac{sin(\theta_j+\theta)}{(1+a^2-2a\cos(\theta_j+\theta)+x)^{(2n+3)/2}}\nonumber\\
& = & \sum_{0\leq j\leq \left[\frac{N-1}{2}\right]} \frac{sin(\theta_j+\theta)}{(1+a^2-2a\cos(\theta_j+\theta)+x)^{(2n+3)/2}}+\sum_{\left[\frac{N-1}{2}\right]+1\leq j\leq N-1} \frac{sin(\theta_j+\theta)}{(1+a^2-2a\cos(\theta_j+\theta)+x)^{(2n+3)/2}}\nonumber\\
& = & \sum_{0\leq j\leq \left[\frac{N-1}{2}\right]} \frac{sin(\theta_j+\theta)}{(1+a^2-2a\cos(\theta_j+\theta)+x)^{(2n+3)/2}}+\sum_{1\leq k\leq N-1-\left[\frac{N-1}{2}\right]} \frac{sin(\theta_{N-k}+\theta)}{(1+a^2-2a\cos(\theta_{N-k}+\theta)+x)^{(2n+3)/2}}\nonumber
\end{eqnarray}
\begin{eqnarray}
& = & \sum_{0\leq j\leq \left[\frac{N-1}{2}\right]} \frac{sin(\theta_j+\theta)}{(1+a^2-2a\cos(\theta_j+\theta)+x)^{(2n+3)/2}}-\sum_{1\leq k\leq N-1-\left[\frac{N-1}{2}\right]} \frac{sin(\theta_k-\theta)}{(1+a^2-2a\cos(\theta_k-\theta)+x)^{(2n+3)/2}}\nonumber\\
& = & \sum_{0\leq j\leq \left[\frac{N-1}{2}\right]} \frac{sin(\theta_j+\theta)}{(1+a^2-2a\cos(\theta_j+\theta)+x)^{(2n+3)/2}}-\sum_{1\leq k\leq \left[\frac{N}{2}\right]} \frac{sin(\theta_{k}-\theta)}{(1+a^2-2a\cos(\theta_{k}-\theta)+x)^{(2n+3)/2}}\nonumber
\end{eqnarray}\\\\
From Lemma 2.9, we only need to prove that:\\\\
i.~~~~~$d_n=\left[\sum_{0\leq j\leq \left[\frac{N-1}{2}\right]} \frac{sin(\theta_j+\theta)}{(1+a^2-2a\cos(\theta_j+\theta)+x)^{(2n+3)/2}}\right]^{\frac{1}{n}}-\left[\sum_{1\leq k\leq \left[\frac{N}{2}\right]} \frac{sin(\theta_{k}-\theta)}{(1+a^2-2a\cos(\theta_{k}-\theta)+x)^{(2n+3)/2}}\right]^{\frac{1}{n}}>0$\\\\
in $\{x:1\geq x\geq0\}$ for some sufficiently large $n \in \mathbb{N}$.\\\\
ii.~~~~~$e_n=x^2d_n>0$ in $\{x:x\geq1\}$ for some sufficiently large $n \in \mathbb{N}$.\\\\
From Lemma 2.8,we know that $$d_n\rightarrow\frac{1}{1+a^2-2a\cos(\theta)+x}-\frac{1}{1+a^2-2a\cos(\frac{2\pi}{N}-\theta)+x},\mbox{when}~ n\rightarrow\infty.$$
But
$$\frac{1}{1+a^2-2a\cos(\theta)+x}-\frac{1}{1+a^2-2a\cos(\frac{2\pi}{N}-\theta)+x}>0,\mbox{when}~\theta \in (0,\frac{\pi}{N})\mbox{in}~\{x:1\geq x\geq0\}$$\\
So $d_n>0$ in $\{x:1\geq x\geq0\}$ for sufficiently large $n\in \mathbb{N}$.\\
Similarly we have\\ $$e_n\rightarrow\frac{x^2}{1+a^2-2a\cos(\theta)+x}-\frac{x^2}{1+a^2-2a\cos(\frac{2\pi}{N}-\theta)+x},\mbox{when}~ n\rightarrow\infty$$\\
However it holds\\
$$\frac{x^2}{1+a^2-2a\cos(\theta)+x}-\frac{x^2}{1+a^2-2a\cos(\frac{2\pi}{N}-\theta)+x}=\frac{2ax^2\left(\cos(\theta)-\cos(\frac{2\pi}{N}-\theta)\right)}{(1+a^2-2a\cos(\theta)+x)(1+a^2-2a\cos(\frac{2\pi}{N}-\theta)+x)}>0,$$\\
when $\theta \in (0,\frac{\pi}{N})$ in $\{x:x\geq1\}$\\
So $e_n>0$ in $\{x:x\geq1\}$ for sufficiently large $n \in \mathbb{N}$.\\
As a result $$g_n(x)=\sum_{0\leq j\leq N-1} \frac{sin(\theta_j+\theta)}{(1+a^2-2a\cos(\theta_j+\theta)+x)^{(2n+3)/2}}>0$$
for sufficiently large $n \in \mathbb{N}$ for any $\theta \in (0,\frac{\pi}{N})$ in $\{x:x\geq0\}$.$~~~~~~~~~~~~~~~~~~~~~~~~~~~~~~~~~~~~~~~~~~~~~~~~~~~~~~~~~~~~~~~~~~~~~~~\Box$\\

\setcounter{equation}{0}
\section{The Proofs of Main Results}
{\bf Proof of Theorem 1.8:}\\\\
1) When $h=0$ and $a=1$,in order to avoid overlap, we assume $\theta\neq0,2\pi/N$.So we only consider $\theta\in (0,2\pi/N)$.From Lemma 2.10, we know that,if we can prove $$f(\theta)=\sum_{1\leq j\leq N} \frac{\sin(\theta_j+\theta)}{|1+a^2-2a\cos(\theta_j+\theta)+h^2|^{3/2}}>0$$ for any $a>0$, $h\geq0$ and $\theta\in (0,\frac{\pi}{N})$,then there must be $\theta=\pi/N$.\\But by Lemma 2.11,it is obvious.\\\\
2)~~~Except for $h=0$ and $a=1$, from Lemma 2.10 and Lemma 2.11,we know that there must be $\theta=0\mbox{ or }\theta=\pi/N.~~~~~~~~~~~~~~~~~~~~~~~~~~~~~~~~~~~~~~~~~~~~~~~~~~~~~~~~~~~~~~~~~~~~~~~~~~~~~~~~~~~~~~~~~~~~~~~~~~~~~~~~~~~~~~~~~~~~~~~~~~~~~~~~~~~~~~~~~~~~~~~\Box$\\\\
{\bf Proof of Corollary 1.9:}\\

From Theorem 1.8 we know $\theta=\pi/N$,thus we have $b=1$ by Corollary1.5. We know it is really a central configuration because of Lemma 2.6.$~~~~~~~~~~~~~~~~~~~~~~~~~~~~~~~~~~~~~~~~~~~~~~~~~~~~~~~~~~~~~~~~~~~~~~~~~~~~~~~~~~~~~~~~~~~~~~~~~~~\Box$\\\\
{\bf Proof of Corollary1.10:}\\

It is obvious from Lemma 2.4 ,Lemma 2.5 and Theorem 1.8. $~~~~~~~~~~~~~~~~~~~~~~~~~~~~~~~~~~~~~~~~~~~~~~~~~~~~~~~~~~~~~~~~~~~~~~~\Box$\\\\
{\bf Proof of Corollary 1.11:}\\

We have the following equations from Lemma 2.5 for $a=1$:\\
\begin{eqnarray}
A-b\sum_{1\leq j\leq N} \frac{\cos(\theta_j+\theta)}{(2-2\cos(\theta_j+\theta)+h^2)^{3/2}}=\frac{\mu}{1+b}\\\nonumber\\
bA-\sum_{1\leq j\leq N} \frac{\cos(\theta_j+\theta)}{(2-2\cos(\theta_j+\theta)+h^2)^{3/2}}=\frac{b\mu}{1+b}\\\nonumber\\
\sum_{1\leq j\leq N} \frac{sin(\theta_j+\theta)}{(2-2\cos(\theta_j+\theta)+h^2)^{3/2}}=0~~~~~~\\\nonumber\\
\sum_{1\leq j\leq N} \frac{1}{(2-2\cos(\theta_j+\theta)+h^2)^{3/2}}=\frac{\mu}{1+b}~~~~~
\end{eqnarray}\\
By (3.1) and (3.2) we have
\begin{eqnarray}
(b^2-1)\sum_{1\leq j\leq N} \frac{\cos(\theta_j+\theta)}{(2-2\cos(\theta_j+\theta)+h^2)^{3/2}}=0\\\nonumber\\
(b^2-1)A=\frac{(b^2-1)\mu}{1+b}~~~~~~~~~~~~~~~~~
\end{eqnarray}\\
i.if $b\neq1$,then\\
\begin{eqnarray}
\sum_{1\leq j\leq N} \frac{\cos(\theta_j+\theta)}{(2-2\cos(\theta_j+\theta)+h^2)^{3/2}}=0\\
A=\frac{\mu}{1+b}~~~~~~~~~~~~~~~~~~~
\end{eqnarray}\\
From (3.3) and (3.7) we can get
\begin{equation}
\sum_{1\leq j\leq N} \frac{\sin(\theta_j)}{(2-2\cos(\theta_j+\theta)+h^2)^{3/2}}=0
\end{equation}
However, we have\\
\begin{flalign}
&~~~~~~~\sum_{1\leq j\leq N} \frac{\sin(\theta_j)}{(2-2\cos(\theta_j+\theta)+h^2)^{3/2}}\nonumber&\\
&~~~~~~~=\sum_{1\leq j\leq \left[\frac{N}{2}\right]} \frac{\sin(\theta_j)}{(2-2\cos(\theta_j+\theta)+h^2)^{3/2}}+\sum_{\left[\frac{N}{2}\right]+1\leq j\leq N} \frac{\sin(\theta_j)}{(2-2\cos(\theta_j+\theta)+h^2)^{3/2}}\nonumber&\\
&~~~~~~~=\sum_{1\leq j\leq \left[\frac{N}{2}\right]} \frac{\sin(\theta_j)}{(2-2\cos(\theta_j+\theta)+h^2)^{3/2}}+\sum_{0\leq k\leq N-1-\left[\frac{N}{2}\right]} \frac{sin(\theta_{N-k})}{(2-2\cos(\theta_{N-k}+\theta)+h^2)^{3/2}}\nonumber&\\
&~~~~~~~=\sum_{1\leq j\leq \left[\frac{N}{2}\right]} \frac{\sin(\theta_j)}{(2-2\cos(\theta_j+\theta)+h^2)^{3/2}}+\sum_{1\leq k\leq \left[\frac{N-1}{2}\right]} \frac{\sin(\theta_{N-k})}{(2-2\cos(\theta_{N-k}+\theta)+h^2)^{3/2}}\nonumber&\\
&~~~~~~~=\sum_{1\leq j\leq \left[\frac{N}{2}\right]} \frac{\sin(\theta_j)}{(2-2\cos(\theta_j-\theta)+h^2)^{3/2}}-\sum_{1\leq k\leq \left[\frac{N-1}{2}\right]} \frac{\sin(\theta_k)}{(2-2\cos(\theta_k-\theta)+h^2)^{3/2}}\nonumber&\\
&~~~~~~~=\sum_{1\leq k\leq \left[\frac{N-1}{2}\right]} \left[\frac{\sin(\theta_k)}{(2-2\cos(\theta_k+\theta)+h^2)^{3/2}}-\frac{\sin(\theta_k)}{(2-2\cos(\theta_k-\theta)+h^2)^{3/2}}\right]&\nonumber\\
&~~~~~~~~~~+\sum_{\left[\frac{N-1}{2}\right]< k\leq \left[\frac{N}{2}\right]} \frac{\sin(\theta_k)}{(2-2\cos(\theta_k-\theta)+h^2)^{3/2}}&
\end{flalign}\\
It is easy to know that no matter what $N$ is even or odd, we can get\\
\begin{equation}\sum_{\left[\frac{N-1}{2}\right]< k\leq \left[\frac{N}{2}\right]} \frac{\sin(\theta_k)}{(2-2\cos(\theta_k-\theta)+h^2)^{3/2}}=0\nonumber
\end{equation}\\
Because when $N$ is even, from $\left[\frac{N-1}{2}\right]< k\leq\left[\frac{N}{2}\right]$,~$k$ has to be $\frac{N}{2}$,then $\sin\theta_k=0$. When $N$ is odd, from $\left[\frac{N-1}{2}\right]< k\leq\left[\frac{N}{2}\right]$,~$k$ does not exist.\\\\
So (3.10) must be\\
\begin{flalign}
&~~~~~~~\sum_{\left[1\leq j\leq N\right]} \frac{\sin(\theta_j)}{(2-2\cos(\theta_j+\theta)+h^2)^{3/2}}&\nonumber\\
	&~~~~~~~=\sum_{1\leq k\leq \left[\frac{N-1}{2}\right]} \left[\frac{\sin(\theta_k)}{(2-2\cos(\theta_k+\theta)+h^2)^{3/2}}-\frac{\sin(\theta_k)}{(2-2\cos(\theta_k-\theta)+h^2)^{3/2}}\right]&
\end{flalign}\\
\begin{flalign}
&\mbox{Let}~f_k(\theta)=\frac{\sin(\theta_k)}{(2-2\cos(\theta_k+\theta)+h^2)^{3/2}}-\frac{\sin(\theta_k)}{(2-2\cos(\theta_k-\theta)+h^2)^{3/2}}&
\end{flalign}
Then
\begin{equation}
f'_k(\theta)=-3\sin(\theta_k)\left[\frac{\sin(\theta_k+\theta)}{(2-2\cos(\theta_k+\theta)+h^2)^{5/2}}+\frac{\sin(\theta_k-\theta)}{(2-2\cos(\theta_k-\theta)+h^2)^{5/2}}\right]
\end{equation}\\
For every $k$ satisfying $1\leq k\leq \left[\frac{N-1}{2}\right]$,it is easy to know that $\theta_k+\theta \in (0,\pi]$ and $\theta_k-\theta \in (0,\pi]$ for every $\theta \in [-\frac{\pi}{N},\frac{\pi}{N}]$.\\
Furthermore, both $\theta_k+\theta$ and $\theta_k-\theta$ can't be $\pi$ at the same time.So $f'_k(\theta)<0$ for every $\theta \in [-\frac{\pi}{N},\frac{\pi}{N}]$.\\
Hence we know that (3.9) has at most one solution. But we know that $\theta=0$ is one solution of (3.9).\\
As a result we have $\theta=0$.\\
Then by (3.4) and (3.7) we have\\
$$\sum_{1\leq j\leq N} \frac{1-\cos(\theta_j)}{(2-2\cos(\theta_j)+h^2)^{3/2}}=\frac{\mu}{1+b}
$$\\
or\\
$$\sum_{1\leq j\leq N} \frac{2\sin^2(\theta_j/2)}{(4\sin^2(\theta_j/2)+h^2)^{3/2}}=\frac{\mu}{1+b}$$\\
So $$\frac{\mu}{1+b}<\sum_{1\leq j\leq N} \frac{2\sin^2(\theta_j/2)}{(4\sin^2(\theta_j/2))^{3/2}}=A,$$
But this contradict (3.8).\\
Hence it must be\\
ii.~~~$b=1$.Then (3.1)-(3.4) turn into
\begin{eqnarray}
A-\sum_{1\leq j\leq N} \frac{\cos(\theta_j+\theta)}{(2-2\cos(\theta_j+\theta)+h^2)^{3/2}}=\frac{\mu}{2}\\
\sum_{1\leq j\leq N} \frac{\sin(\theta_j+\theta)}{(2-2\cos(\theta_j+\theta)+h^2)^{3/2}}=0~~~\\
\sum_{1\leq j\leq N} \frac{1}{(2-2\cos(\theta_j+\theta)+h^2)^{3/2}}=\frac{\mu}{2}~~
\end{eqnarray}\\
Thus we have\\
\begin{equation}
\sum_{1\leq j\leq N} \frac{\cos(\theta_j+\theta)}{(2-2\cos(\theta_j+\theta)+h^2)^{3/2}}+\sum_{1\leq j\leq N} \frac{1}{(2-2\cos(\theta_j+\theta)+h^2)^{3/2}}-A=0
\end{equation}
Let
\begin{equation}
g(h)=\sum_{1\leq j\leq N} \frac{\cos(\theta_j+\theta)}{(2-2\cos(\theta_j+\theta)+h^2)^{3/2}}+\sum_{1\leq j\leq N} \frac{1}{(2-2\cos(\theta_j+\theta)+h^2)^{3/2}}-A
\end{equation}
Then
\begin{eqnarray}
g'_n(h) & = & \sum_{1\leq j\leq N} \frac{-3h\cos(\theta_j+\theta)}{(2-2\cos(\theta_j+\theta)+h^2)^{5/2}}+\sum_{1\leq j\leq N} \frac{-3h}{(2-2\cos(\theta_j+\theta)+h^2)^{3/2}}\nonumber\\
& = & -3h\sum_{1\leq j\leq N} \frac{\cos(\theta_j+\theta)+1}{(2-2\cos(\theta_j+\theta)+h^2)^{5/2}}\nonumber
\end{eqnarray}
It is easy to know that $g'(h)<0$ for $h>0$ and every $\theta$.\\
However we have\\
\begin{equation}
\lim_{h\rightarrow\infty} g(h)=\lim_{h\rightarrow\infty} \left(\sum_{1\leq j\leq N} \frac{\cos(\theta_j+\theta)+1}{(2-2\cos(\theta_j+\theta)+h^2)^{3/2}}-A\right)=-A<0
\end{equation}\\
When $\theta=0$ or $\theta=2\pi/N$, it is easy to know that\\
\begin{equation}
\lim_{h\rightarrow0} g(h)=\lim_{h\rightarrow0} \left(\sum_{1\leq j\leq N} \frac{\cos(\theta_j+\theta)+1}{(2-2\cos(\theta_j+\theta)+h^2)^{3/2}}-A\right)=\infty\\
\end{equation}\\
When $\theta=\pi/N$, from (1.27) we have\\
\begin{equation}
 \lim_{h\rightarrow0} g(h) = \sum_{1\leq j\leq N} \frac{\cos(\theta_j+\frac{\pi}{N})+1}{(2-2\cos(\theta_j+\frac{\pi}{N}))^{3/2}}-A>0
\end{equation}\\

So there is exactly one $h>0$,which satisfies (3.17) for $\theta=0$ or $\theta=\pi/N$.

In conclusion, if $a=1$,then $b=1$,$\theta=0$ or $\theta=\pi/N$,there exists a unique $h$ for $\theta=0$ or $\theta=\pi/N$.In other words, there are exactly two spatial central configurations of parallel regular N-gons with the same sizes.$~~~~~~~~~~~~~~~~~~~\Box$

\end{document}